\documentstyle[12pt]{article}

\newcommand{\const}{\mathop{\rm const}\limits}

\newcommand{\supp}{\mathop{\rm supp}\limits}

\begin{document}

\begin{center}

{\bf MULTIPLE WEIGHT RIESZ AND FOURIER TRANSFORMS }

\vspace{3mm}

{\bf IN BILATERAL ANISOTROPIC GRAND LEBESGUE SPACES} \\

\vspace{3mm}

 $ {\bf E.Ostrovsky^a, \ \ L.Sirota^b } $ \\

\vspace{4mm}

$ ^a $ Corresponding Author. Department of Mathematics and computer science, Bar-Ilan University, 84105, Ramat Gan, Israel.\\
\end{center}
E - mail: \ eugostrovsky@list.ru; \ galo@list.ru \\
\begin{center}
$ ^b $  Department of Mathematics and computer science. Bar-Ilan University,
84105, Ramat Gan, Israel.\\
\end{center}
E - mail: \ sirota@zahav.net.il\\

\vspace{3mm}

{\bf Abstract}.  In this short article we introduce so-called
anisotropic (weight) Grand Lebesgue Spaces (more exactly,
Grand Lebesgue-Riesz Spaces), which are generalization of the
classical Lebesgue-Riesz Spaces and ordinary Grand Lebesgue Spaces,
and investigate the boundedness  of weight Riesz potential and
weight Fourier transforms in this spaces.\par

 We construct also several examples to show the exactness of offered estimations.\par

\vspace{3mm}

\par

\vspace{3mm}

\vspace{3mm}

\section{Introduction. Notations. Problem Statement.} \par
\vspace{3mm}

   We will describe now using rearrangement invariant Banach spaces
 of measurable functions.\par

\vspace{3mm}

 {\bf 1. Classical Lebesgue-Riesz Spaces.} \par

 \vspace{3mm}

  Let $ (X,A,\mu) $ be measurable space with sigma-finite
non - trivial measure $ \mu. $ For the measurable real valued
function $ f(x), \ x \in X, f: X \to R $ the symbol
$ |f|_p = |f|_p(X,\mu) $ will denote the usually $ L_p $ norm:
$$
|f|_p = |f|L_p(X, \mu) = \left[ \int_X |f(x)|^p \ \mu(dx) \right]^{1/p}, \ p \ge 1;
$$

 $$
  L_p(X,\mu) = L_p = \{ f, \ f: X \to R, |f|_p < \infty  \}.
 $$ \par

 We will use further only the case when $ X $  is classical  whole Euclidean space $ X = R^d  $
or some its measurable subset with positive measure and $ \mu $ is ordinary Lebesgue measure.  \par

\vspace{3mm}

{\bf 2. Grand Lebesgue Spaces (GLS).}\par

\vspace{3mm}

 We recall in this section  for reader conventions some definitions and facts from the theory
of GLS spaces.\par

\vspace{3mm}

 Recently, see \cite{Fiorenza1}, \cite{Fiorenza2}, \cite{Fiorenza3}, \cite{Iwaniec1}, \cite{Iwaniec2},
 \cite{Kozachenko1},\cite{Liflyand1}, \cite{Ostrovsky1}, \cite{Ostrovsky2},   etc.
 appears the so-called Grand Lebesgue Spaces $ GLS = G(\psi) =G\psi =
 G(\psi; A,B), \ A,B = \const, A \ge 1, A < B \le \infty, $ spaces consisting
 on all the measurable functions $ f: X \to R $ with finite norms

     $$
     ||f||G(\psi) \stackrel{def}{=} \sup_{p \in (A,B)} \left[ |f|_p /\psi(p) \right].
     $$

      Here $ \psi(\cdot) $ is some continuous positive on the {\it open} interval
    $ (A,B) $ function such that

     $$
     \inf_{p \in (A,B)} \psi(p) > 0, \ \psi(p) = \infty, \ p \notin (A,B).
     $$
We will denote
$$
 \supp (\psi) \stackrel{def}{=} (A,B) = \{p: \psi(p) < \infty, \}
$$

The set of all $ \psi $  functions with support $ \supp (\psi)= (A,B) $ will be
denoted by $ \Psi(A,B). $ \par
  This spaces are rearrangement invariant, see \cite{Bennet1}, and
  are used, for example, in the theory of probability  \cite{Kozachenko1},
  \cite{Ostrovsky1}, \cite{Ostrovsky2}; theory of Partial Differential Equations \cite{Fiorenza2},
  \cite{Iwaniec2};  functional analysis \cite{Fiorenza3}, \cite{Iwaniec1},  \cite{Liflyand1},
  \cite{Ostrovsky2}; theory of Fourier series \cite{Ostrovsky1},
  theory of martingales \cite{Ostrovsky2},mathematical statistics  \cite{Sirota2},
 \cite{Sirota4}, \cite{Sirota5}, \cite{Sirota6}, \cite{Sirota7}, \cite{Sirota8};
   theory of approximation \cite{Ostrovsky7}   etc.\par

 Notice that in the case when $ \psi(\cdot) \in \Psi(A,\infty)  $ and a function
 $ p \to p \cdot \log \psi(p) $ is convex,  then the space
$ G\psi $ coincides with some {\it exponential} Orlicz space. \par
 Conversely, if $ B < \infty, $ then the space $ G\psi(A,B) $ does  not coincides with
 the classical rearrangement invariant spaces: Orlicz, Lorentz, Marcinkiewicz  etc.\par

\vspace{3mm}

{\bf Remark 1.1} If we introduce the {\it discontinuous} function

$$
\psi_r(p) = 1, \ p = r; \psi_r(p) = \infty, \ p \ne r, \ p,r \in (A,B)
$$
and define formally  $ C/\infty = 0, \ C = \const \in R^1, $ then  the norm
in the space $ G(\psi_r) $ coincides with the $ L_r $ norm:

$$
||f||G(\psi_r) = |f|_r.
$$
Thus, the Grand Lebesgue Spaces are direct generalization of the
classical exponential Orlicz's spaces and Lebesgue spaces $ L_r. $ \par

\vspace{3mm}

{\bf Remark 1.2}  The function $ \psi(\cdot) $ may be generated as follows. Let $ \xi = \xi(x)$
be some measurable function: $ \xi: X \to R $ such that $ \exists  (A,B):
1 \le A < B \le \infty, \ \forall p \in (A,B) \ |\xi|_p < \infty. $ Then we can
choose

$$
\psi(p) = \psi_{\xi}(p) = |\xi|_p.
$$

 Analogously let $ \xi(t,\cdot) = \xi(t,x), t \in T, \ T $ is arbitrary set,
be some {\it family } $ F = \{ \xi(t, \cdot) \} $ of the measurable functions:
$ \forall t \in T  \ \xi(t,\cdot): X \to R $ such that
$$
 \exists  (A,B): 1 \le A < B \le \infty, \ \sup_{t \in T} \
|\xi(t, \cdot)|_p < \infty.
$$
Then we can choose

$$
\psi(p) = \psi_{F}(p) = \sup_{t \in T}|\xi(t,\cdot)|_p.
$$
The function $ \psi_F(p) $ may be called as a {\it natural function} for the family $ F. $
This method was used in the probability theory, more exactly, in
the theory of random fields, see \cite{Ostrovsky1}. \par

\vspace{3mm}

{\bf 3. Anisotropic Lebesgue-Riesz spaces.}

\vspace{3mm}

 We recall here the definition of
 the so-called anisotropic Lebesgue (Lebesgue-Riesz) spaces. More detail information about this
spaces see in the books  of Besov O.V., Il'in V.P., Nikol'skii S.M.
\cite{Besov1}, chapter 16,17; Leoni G. \cite{Leoni1}, chapter 11; using for us theory of
operators interpolation in this spaces see in \cite{Besov1}, chapter 17,18. \par

  Let $ (X_j,A_j,\mu_j) $ be measurable spaces with sigma-finite
non - trivial measures $ \mu_j. $
Let also $ p = (p_1, p_2, . . . , p_l) $ be $ l- $ dimensional vector such that
$ 1 \le p_j \le \infty.$ \par

 Recall that the anisotropic Lebesgue space $ L_{ \vec{p}} $ consists on all the  total measurable
real valued function  $ f = f(x_1,x_2,\ldots, x_l) = f( \vec{x} ) $

$$
f:  \otimes_{j=1}^l X_j \to R
$$

with finite norm $ |f|_{ \vec{p} } \stackrel{def}{=} $

$$
\left( \int_{X_l} \mu_l(dx_l) \left( \int_{X_{l-1}} \mu_{l-1}(dx_{l-1}) \ldots \left( \int_{X_1}
 |f(\vec{x})|^{p_1} \mu(dx_1) \right)^{p_2/p_1 }  \ \right)^{p_3/p_2} \ldots   \right)^{1/p_l}.
$$

 Note that in general case $ |f|_{p_1,p_2} \ne |f|_{p_2,p_1}, $
but $ |f|_{p,p} = |f|_p. $ \par

 Observe also that if $ f(x_1, x_2) = g_1(x_1) \cdot g_2(x_2) $ (condition of factorization), then
$ |f|_{p_1,p_2} = |g_1|_{p_1} \cdot |g_2|_{p_2}, $ (formula of factorization). \par

\vspace{3mm}

{\bf 4. Anisotropic Grand Lebesgue-Riesz spaces.}

\vspace{3mm}
 Let $ Q $ be convex (bounded or not) subset of the set $ \otimes_{j=1}^l [1,\infty]. $
 Let $ \psi = \psi(\vec{p}) $ be continuous in an interior $ Q^0 $ of the set $ Q $
strictly  positive  function such that

$$
\inf_{\vec{p} \in Q^0}  \psi(\vec{p}) > 0; \ \inf_{\vec{p} \notin Q^0}  \psi(\vec{p}) = \infty.
$$

 We denote the set all of such a functions as $ \Psi(Q). $ \par
The  Anisotropic Grand Lebesgue Spaces $ AGLS = AGLS(\psi) $
 space consists on all the measurable functions

$$
f:  \otimes_{j=1}^l X_j \to R
$$
with finite (mixed) norms

$$
||f||AG\psi = \sup_{\vec{p} \in Q^0} \left[ \frac{|f|_{\vec{p}}}{\psi(\vec{p} )} \right].
$$

 An application  into the theory of multiple Fourier transform of these  spaces see in articles \cite{Benedek1} and
 \cite{Ostrovsky3}, where are considered some problems
of boundedness of singular operators  in (weight) Grand Lebesgue Spaces and  Anisotropic Grand Lebesgue Spaces.
We intend to generalize some  results obtained in \cite{Benedek1}, \cite{Ostrovsky3}. \par

\vspace{3mm}

{\bf 5. Description of considered operators.}

\vspace{3mm}

{\bf 5A. Weight Riesz's ordinary operator.} Recall that the operator, more exactly,
the {\it family } of operators $ I_{\alpha,\beta,\gamma}[f](x), \ x \in R^d $ of a view

 $$
u(x) = I_{\alpha,\beta,\gamma}[f](x) =I_{d;\alpha,\beta,\gamma}[f](x) = I[f](x) = |x|^{-\beta}
\int_{R^d} \frac{f(y) \ |y|^{-\alpha} \ dy}{|x - y|^{\gamma}};  \eqno(1.5A.1)
 $$
is said to be {\it weight} Riesz's operator;
here $ \gamma = \const, \alpha,\beta = \const \ge 0, \alpha + \gamma < d  $
and $ |x| $ denotes ordinary Euclidean norm of the vector $ x. $ \par

 We refer used in this article  results from \cite{Ostrovsky4}.

 Let us define the function $ q=q(p) $ as follows:
$$
d+\frac{d}{q} = \frac{d}{p} + (\alpha+\beta+\gamma). \eqno(1.5A.2)
$$

 We will denote the set of all such a values $ (p,q) $  as $ G(d;\alpha,\beta,\gamma) $ or
for simplicity $ G = G(d;\alpha,\beta,\gamma). $  \par
 Further we will suppose in this subsection that
$ (p,q) \in G(d;\alpha,\beta,\gamma) = G. $ \par

  We denote also
$$
p_-:=\frac{d}{d-\alpha}, \ p_+:= \frac{d}{d-\alpha-\gamma};
$$
and correspondingly
$$
q_-:= \frac{d}{\beta+\gamma}, \ q_+ := \frac{d}{\beta},
$$
where in the case $ \beta=0 \ \Rightarrow q_+:= + \infty; $
$$
\kappa \stackrel{def}{=} \kappa(\alpha,\beta,\gamma) := (\alpha + \beta + \gamma)/d.
$$

 There exist a positive finite coefficient $  K_{R;d; \alpha,\beta, \gamma}(p), \ p \in (p_-,p_+)  $ for which


$$
|I_{\alpha,\beta,\gamma}[f]|_q  \le K_{R;d; \alpha,\beta, \gamma}(p), \eqno(1.5A.3)
$$
where by definition

$$
K_{R;d; \alpha,\beta, \gamma}(p) =  \sup_{ p \in (p_-, p_+) } \sup_{0 < |f|_p   < \infty}
\left[\frac{|I_{\alpha,\beta,\gamma}[f]|_{q(p)}}{ |f|_p } \right]. \eqno(1.5A.4)
$$
 It is known that the coefficient $ K_{R;d; \alpha,\beta, \gamma}(p) $ is bilateral bounded :

 $$
\frac{C_1(d; \alpha,\beta, \gamma)}{ \left[(p-p_-)(p_+ - p) \right]^{\kappa}}  \le K_R(d; \alpha,\beta, \gamma) \le 
\frac{C_2(d; \alpha,\beta, \gamma)} {\left[(p-p_-)(p_+ - p) \right]^{\kappa} }. \eqno(1.5A.5)
 $$

\vspace{3mm}

 Notice  that this estimation may be obtained also by  the direct computation from the recent article
of Cruz-Uribe D.,  Moen K. \cite{Cruz-Uribe1}.

\vspace{3mm}

{\bf 5B. Weight Fourier transform.}\\

\vspace{3mm}

 We define alike to the book Okikiolu  \cite{Okikiolu1}, p. 313-314, see also \cite{Beckner2}
 the so-called {\it double weight} Fourier transform $ F_{\alpha,\beta}[f](x), \ x \in R^d $  by the following way:

$$
F_{\alpha,\beta}[f](x) = (2 \pi)^{-d/2} |x|^{\alpha} \int_{R^d} |y|^{\beta} \ f(y) \ e^{i x y} \ dy, \eqno(1.5B.1)
 $$

$$
F[f](x) = F_{0,0}[f](x) = (2 \pi)^{-d/2} \int_{R^d}  \ f(y) \ e^{i x y} \ dy,
 $$
where $ xy $ denotes the inner product of the vectors $ x, y : xy = \sum_k x_k y_k.  $ \par

 Let us introduce the following important conditions in the domain of parameters $ (p,q): $

$$
p_0: = d/(d - \beta), \ p_1: = \infty, \ q_0 := 1, \  q_1: = d/\alpha >1;
\beta - \alpha = d - d(1/p + 1/q), \eqno(1.5B.2)
$$

$$
p > p_0, \ 1 \le q < q_1. \eqno(1.5B.3)
$$
which defined uniquely the continuous function $ q = q(p). $ \par

 The inequality of a view (relative the Lebesgue measures in whole space $ R^d $ )

$$
|F_{\alpha,\beta}[f]|_q \le K_{PBO}(p) \ |f|_p \eqno(1.5B.4)
$$
is said to be (generalized) weight Pitt - Beckner - Okikiolu (PBO) inequality.\par

 We will understood as customary as the value $ K_{PBO}(p) $ its minimal value:

$$
K_{PBO}(p) = \sup_{f \ne 0, |f|_p < \infty} \left[ \frac{|F_{\alpha,\beta}[f]|_q}{|f|_p} \right].
\eqno(1.5B.5)
$$

 It is known \cite{Ostrovsky3} that

$$
C_1(\alpha,\beta,d)  \left[\frac{p}{p-p_0} \right]^{(\alpha +\beta)/d } \le K_{PBO}(p) \le
C_2(\alpha,\beta,d)  \left[\frac{p}{p-p_0} \right]^{ \max(1,(\alpha +\beta)/d) }, \eqno(1.5B.6)
$$
$ 0 < C_{1,2}(\alpha,\beta,d) < \infty.  $

 The right hand side of bilateral estimate (1.5B.6) may be obtained after  some computations
 from the articles \cite{Benedetto1}, \cite{Benedetto2}, \cite{Beckner2}, \cite{Carton-Lebrun1},
\cite{Heining1}, \cite{Jurkat1}.  The lower estimate is obtained in \cite{Ostrovsky3}. \par

 See also many publications about this problem
\cite{Lacey1}, \cite{Lakey1}, \cite{Lakey2},\cite{Muckenhoupt1}, \cite{Muckenhoupt2},
\cite{P'erez1}, \cite{Petermichl1}, \cite{Samko1}, \cite{Schechter1}, \cite{Stein1},
\cite{Tomas-Stein1} etc. \par

\vspace{4mm}

 The paper is organized as follows. In the next section we investigate multiple weight Riesz and
Fourier transforms  and ground the neediness of introducing of anisotropic spaces.
 In the third section we obtain the conditions for boundedness  of the
multiple weight Riesz and Fourier transforms in anisotropic Grand Lebesgue-Riesz spaces.
 The fourth section is devoted to the  multiple Riesz potential in bounded domains.
  The last section contains some slight generalizations and concluding remarks.

\vspace{4mm}

 We use symbols $C(X,Y),$ $C(p,q;\psi),$ etc., to denote positive finite
constants along with parameters they depend on, or at least
dependence on which is essential in our study. To distinguish
between two different constants depending on the same parameters
we will additionally enumerate them, like $C_1(X,Y)$ and $C_2(X,Y).$

\hfill $\Box$ \\
 \bigskip

\section{ Multiple weight Riesz and Fourier transforms. Neediness of anisotropic spaces.}

\vspace{3mm}

 We consider in this section the weight multidimensional (vector): $ l \ge 2 $ generalization of
weight Riesz's $ L_p $ estimation and PBO inequality.\par
 In this section $ x = \vec{x} \in R^d $ be $ d- $ dimensional vector, $ d = 1, 2, . . .  $ which consists
on the $ l $ subvectors  $ {x_j}, j = 1, 2, . . . , l : \ $
$$
x = (x_1, x_2, . . . , x_l), \ x_j = \vec{x_j} \in  R^{m_j} , \dim(x_j) = m_j \ge 1;
$$
$$
 \alpha = \vec{\alpha} = \{\alpha_1, \alpha_2, \ldots, \alpha_l \},
 \beta = \vec{\beta} = \{\beta_1,\beta_2, \ldots,\beta_l \}
$$
  and   $ \gamma = \vec{\gamma} = \{\gamma_1,\gamma_2, \ldots,\gamma_l \}  $
  be three fixed $ l - $ dimensional numerical vectors. \par
   Note that in the case of weight multiple Riesz potential the third vector
$\vec{\gamma} $ is absent. \par

 We denote as ordinary
$$
|x|^{-\alpha} = |\vec{x}|^{-\vec{\alpha}} = \prod_{j=1}^l x_j^{-\alpha_j}, \
|y|^{\beta} = |\vec{y}|^{\vec{\beta}} = \prod_{j=1}^l y_j^{\beta_j}, \
|z|^{\gamma} = |\vec{z}|^{\vec{\gamma}} = \prod_{j=1}^l z_j^{\gamma_j}. \eqno(2.1.0)
$$
 Obviously, $  \sum_j m_j = d. $\par

\vspace{3mm}

{\bf A. Weight multiple Riesz potential.}

\vspace{3mm}

 Let $ f, \ f: R^d \to R $ be (total) measurable function.
 Let $ \gamma_j = \const, \alpha_j,\beta_j = \const \ge 0, \alpha_j + \gamma_j < m_j, \ j=1,2,\ldots,l.  $

We define the function $ q_j=q_j(p_j) $ as follows:
$$
m_j \left(1+\frac{1}{q_j} - \frac{1}{p_j} \right) = \alpha_j+\beta_j+\gamma_j, \eqno(2.1.1)
$$
i.e. $ (p_j,q_j) \in G(m_j; \alpha_j,\beta_j,\gamma_j) =: G_j. $ \par

 Further we will suppose in this subsection that
$ (p_j,q_j) \in G(m_j;\alpha_j,\beta_j,\gamma_j) =: G_j. $ \par

  We denote also as before
$$
p_-^{(j)}:=\frac{m_j}{m_j-\alpha_j}, \ p_+^{(j)}:= \frac{m_j}{m_j-\alpha_j-\gamma_j}; \eqno(2.1.2)
$$
and correspondingly
$$
q_-^{(j)}:= \frac{m_j}{\beta_j+\gamma_j}, \ q_+^{(j)} := \frac{m_j}{\beta_j}, \eqno(2.1.2a)
$$
where in the case $ \beta_j=0 \ \Rightarrow q_+^{(j)}:= + \infty; $
$$
\kappa_j \stackrel{def}{=} \kappa_j(\alpha_j,\beta_j,\gamma_j) := (\alpha_j + \beta_j + \gamma_j)/m_j. \eqno(2.1.3)
$$

 We define the following {\it multiple} weight Riesz (linear) potential:

 $$
u(\vec{x}) = I_{\vec{\alpha}, \vec{\beta}, \vec{\gamma}}[f](\vec{x})  =
 I_{\vec{m};\vec{\alpha}, \vec{\beta}, \vec{\gamma}}[f](\vec{x}) =
\otimes_{j=1}^l I_{m_j; \alpha_j,\beta_j,\gamma_j}[f](\vec{x})  =
$$

$$
 \int_{R^{m_1}}  \ \frac{|x_1|^{-\beta_1} \ |y_1|^{-\alpha_1} \ dy_1}{|x_1 - y_1|^{\gamma_1}}
\left[  \int_{R^{m_2}}  \ \frac{|x_2|^{-\beta_2} \ |y_2|^{-\alpha_2} \ dy_2}{|x_2 - y_2|^{\gamma_2}} \left[ \ldots
\left[ \int_{R^{m_l}}  \ \frac{|x_l|^{-\beta_l} \ |y_l|^{-\alpha_l} \ f( \vec{y} ) \ dy_l}{|x_l - y_l|^{\gamma_l}} \right]
 \right]   \right]. \eqno(2.1.4)
$$

  Note that our definition different on the ones in the articles \cite{Heining2}, \cite{Lacey1}, \cite{Moen1},
\cite{Pradolini1}  etc. \par
 More general case, namely when the weight function is  regular varying will be considered further. \par

\vspace{3mm}

{\bf Theorem 2.1.} \par
 The conditions (2.1.1) are necessary and sufficient for the existence non-trivial coefficient
$ K_{{\vec{m}; \vec{\alpha}, \vec{\beta}, \vec{\gamma}}}(\vec{p}) $ for the following estimate:

$$
| I_{\vec{m};\vec{\alpha}, \vec{\beta}, \vec{\gamma}}[f](\vec{\cdot})|_{\vec{q}} \le
K_{{\vec{m}; \vec{\alpha}, \vec{\beta}, \vec{\gamma}}}(\vec{p}) \ |f(\cdot)|_{ \vec{p}}, \eqno(2.1.5)
$$
and under this conditions  for the minimal value of coefficient
$  K_{{\vec{m}; \vec{\alpha}, \vec{\beta}, \vec{\gamma}}}(\vec{p}) $ holds the
following inequality:

$$
 \frac{C_3(\vec{m}, \vec{\alpha}, \vec{\beta}, \vec{\gamma})}
 { \prod_{j=1}^l \left[(p_+^{(j)} -p_j)(p_j- p_-^{(j)})  \right]^{\kappa_j} } \le
    K_{{\vec{m}; \vec{\alpha}, \vec{\beta}, \vec{\gamma}}}(\vec{p}) \le
 \frac{C_4(\vec{m}, \vec{\alpha}, \vec{\beta}, \vec{\gamma})}
 { \prod_{j=1}^l \left[(p_+^{(j)} -p_j)(p_j- p_-^{(j)})  \right]^{\kappa_j} }, \eqno(2.1.6)
$$
if for all the values $ j \ \Rightarrow  p_j \in \left( p_-^{(j)}, p_+^{(j)} \right) $
 and $  K_{{\vec{m}; \vec{\alpha}, \vec{\beta}, \vec{\gamma}}}(\vec{p}) = \infty  $ in other case. \par

\vspace{3mm}

  Recall that the  values $ \{ p_j \} $ are variable, in contradiction to the values
$  \{ p_+^{(j)},  p_-^{(j)} \}.  $ \par

\vspace{3mm}

{\bf Proof.}  We can and will suppose without loss of generality $ f \in S(R^d), \ f \ne 0, $ where $ S(R^d) $
denotes Schwartz class of functions in whole space $ R^d. $  We use the {\it multidimensional} generalization of the
so-called {\it dilation } method, see  \cite{Talenti1}, \cite{Stein1}, chapter 3. \par
 Indeed, let for simplicity $ l = 2. $ Let $  f \in S(R^d), \ d = m_1 + m_2, \ f \ne 0 $ be any function for which the inequality
 (2.1.5) there holds. Let $ \lambda_1, \lambda_2 = \const > 0 $ be arbitrary (independent) numbers.  A multidimensional
 dilation operator $ T_{\lambda_1, \lambda_2}[f] $ may be defined as follows:

 $$
 T_{\vec{\lambda}}[f](x_1, x_2) =  T_{\lambda_1, \lambda_2}[f](x_1, x_2) = f(\lambda_1 x_1, \lambda_2 x_2).
 $$
  Evidently  $ T_{\lambda_1, \lambda_2}[f] \in  S(R^l).  $ \par
  We have consequently:

  $$
  | \ T_{\lambda_1, \lambda_2}[f] \ |_{\vec{p}} = \lambda_1^{-m_1/p_1} \ \lambda_2^{-m_2/p_2} \ |f|_{\vec{p}};
  $$

 $$
 I_{\vec{\alpha}, \vec{\beta}, \vec{\gamma} }[T_{\vec{\lambda}}[f]](x_1/\lambda_1, x_2/\lambda_2) =
 \lambda_1^{\alpha_1 + \beta_1 + \gamma_1-m_1 } \  \lambda_2^{\alpha_2 + \beta_2 + \gamma_2-m_2 } \
T_{\vec{\lambda}}[ I_{\vec{\alpha}, \vec{\beta}, \vec{\gamma} }[f]](x_1, x_2);
 $$

  $$
 I_{\vec{\alpha}, \vec{\beta}, \vec{\gamma} }[T_{\vec{\lambda}}[f]](x_1, x_2) =
 \lambda_1^{\alpha_1 + \beta_1 + \gamma_1-m_1 } \  \lambda_2^{\alpha_2 + \beta_2 + \gamma_2-m_2 } \
T_{\vec{\lambda}}[ I_{\vec{\alpha}, \vec{\beta}, \vec{\gamma} }[f]](x_1 \ \lambda_1, x_2 \ \lambda_2);
 $$

 $$
 | I_{\vec{\alpha}, \vec{\beta}, \vec{\gamma} }[T_{\vec{\lambda}}[f](x_1, x_2)] \ |_{\vec{q}} =
 \lambda_1^{\alpha_1 + \beta_1 + \gamma_1-m_1 - m_1/q_1 } \  \lambda_2^{\alpha_2 + \beta_2 + \gamma_2-m_2 - m_2/q_2 } \
 | I_{\vec{\alpha}, \vec{\beta}, \vec{\gamma} }[f] \ |_{\vec{q}}.
 $$

  Substituting into (2.1.5) we obtain

 $$
C_1 \ \lambda_1^{\alpha_1 + \beta_1 + \gamma_1-m_1 - m_1/q_1 } \  \lambda_2^{\alpha_2 + \beta_2 + \gamma_2-m_2 - m_2/q_2 }   \le
C_2 \ \lambda_1^{-m_1/p_1} \ \lambda_2^{-m_2/p_2}.
 $$
  Since the values $ \lambda_1,\lambda_2 $ are arbitrary positive, we obtain the  equalities
(2.1.1). We refer also the reader to the article \cite{Ostrovsky3} for details. \par

 It remains to obtain the estimations  (2.1.6). First of all we will obtain the upper estimate.\par

 It is sufficient again to consider  only the "two-dimensional" case $  l = 2. $  Namely,  define
the function $ u = u(x_1, x_2), \ x_1 \in R^{m_1}, x_2 \in R^{m_2}  $ as follows:
$$
u(x_1, x_2) = \otimes_{j=1}^2 I_{m_j; \alpha_j,\beta_j,\gamma_j}[f](\vec{x})=
$$

$$
 \int_{R^{m_1}}  \ \frac{|x_1|^{-\beta_1} \ |y_1|^{-\alpha_1} \ dy_1}{|x_1 - y_1|^{\gamma_1}} \
\left[ \int_{R^{m_2}}  \ \frac{|x_2|^{-\beta_2} \ |y_2|^{-\alpha_2} \ f(y_1,y_2) \ dy_2}{|x_2 - y_2|^{\gamma_2}} \right]=
$$

$$
 I_{m_1;\alpha_1,\beta_2,\gamma_1}[g](x_1,x_2),
$$
where

$$
g(x_1,x_2)=  \int_{R^{m_2}}  \ \frac{|x_2|^{-\beta_2} \ |y_2|^{-\alpha_2} \ f(y_1,y_2) \ dy_2}{|x_2 - y_2|^{\gamma_2}}.
$$

 We have taking the $ L(q_1) $ norm on the variable $  x_1 $ and using the one-dimensional version for
Riesz's potential

$$
|u(\cdot, x_2)|L_{q_1, R^m_1} \le
\frac{ C_2(m_1; \alpha_1,\beta_1, \gamma_1)}{ [(p_1-p_-^{(1)})(p_+^{(1)} - p_1) ]^{\kappa_1}} \
|g(\cdot,\cdot)| L_{p_1, R^m_1}.
$$

 We use the triangle inequality for the $ L_{p_1}  $ norm, denoting $ h(y_2) = |f(y_1,y_2)|_{p_1,R^d_1}: $

$$
|g(\cdot,\cdot)| L_{p_1, R^m_1} \le
\int_{R^{m_2}}  \ \frac{|x_2|^{-\beta_2} \ |y_2|^{-\alpha_2} \ |f(y_1,y_2)|_{p_1,R^m_1} \ dy_2}{|x_2 - y_2|^{\gamma_2}}=
$$

$$
\int_{R^{m_2}}  \ \frac{|x_2|^{-\beta_2} \ |y_2|^{-\alpha_2} \ h(y_2) \ dy_2}{|x_2 - y_2|^{\gamma_2}}=
I_{m_2; \alpha_2, \beta_2; \gamma_2}[h](x_2),
$$
therefore

$$
|u(\cdot, x_2)|L_{q_1, R^m_1} \le \frac{ C_2(m_1; \alpha_1,\beta_1, \gamma_1)}{ [(p_1-p_-^{(1)})(p_+^{(1)} - p_1) ]^{\kappa_1}} \times
I_{m_2; \alpha_2, \beta_2; \gamma_2}[h](x_2).
$$

 We have taking the $ L(q_2) $ norm on the variable $  x_2 $ and using again the one-dimensional version for
Riesz's potential

$$
|u(\cdot, \cdot)|L_{q_1, R^m_1} L_{q_2, R^m_2} \le \frac{ C_2(m_1; \alpha_1,\beta_1, \gamma_1)}{ [(p_1-p_-^{(1)})(p_+^{(1)} - p_1) ]^{\kappa_1}} \times
\frac{ C_2(m_2; \alpha_2,\beta_2, \gamma_2)}{ [(p_2-p_-^{(2)})(p_+^{(2)} - p_2) ]^{\kappa_2}} \ |h|_{p_2, R^m_2}
$$
or equally

$$
|u(\cdot)|_{\vec{q}} \le
\frac{\prod_{j=1}^l C_2(m_j; \alpha_j,\beta_j, \gamma_j)} {\prod_{j=1}^l [(p_+^{(j)} -p_j)(p_j- p_-^{(j)}) ]^{\kappa_j} } \cdot
|f|_{ \vec{p} }, \eqno(2.1.7)
$$
which is equivalent to the right-hand side inequality of theorem 2.1 with

$$
C_4(\vec{m}, \vec{\alpha}, \vec{\beta}, \vec{\gamma}) = \prod_{j=1}^l C_2(m_j; \alpha_j,\beta_j, \gamma_j). \eqno(2.1.8)
$$
 The {\it lower} estimation for $  | I_{\vec{m};\vec{\alpha}, \vec{\beta}, \vec{\gamma}}[f](\vec{\cdot})|_{\vec{q}} $
may be obtained by means of consideration of factorized  function of a view

$$
f^{(0)}(\vec{x}) = \prod_{j=1}^l f^{(0)}_{m_j}(x_j),  \eqno(2.1.9)
$$
and analogously

$$
C_3(\vec{m}, \vec{\alpha}, \vec{\beta}, \vec{\gamma}) = \prod_{j=1}^l C_1(m_j; \alpha_j,\beta_j, \gamma_j). \eqno(2.1.10)
$$

 More exactly, the equalities (2.1.8) and (2.1.10) imply that the constants $ C_3(\vec{m}, \vec{\alpha}, \vec{\beta}, \vec{\gamma}) $ and
$ C_4(\vec{m}, \vec{\alpha}, \vec{\beta}, \vec{\gamma}) $ in (2.1.6) may be {\it estimates } as follows:

$$
C_3(\vec{m}, \vec{\alpha}, \vec{\beta}, \vec{\gamma}) \ge \prod_{j=1}^l C_1(m_j; \alpha_j,\beta_j, \gamma_j), \eqno(2.1.11)
$$

$$
C_4(\vec{m}, \vec{\alpha}, \vec{\beta}, \vec{\gamma}) \le \prod_{j=1}^l C_2(m_j; \alpha_j,\beta_j, \gamma_j). \eqno(2.1.12)
$$
 The necessity of relations (2.1.1) provided by means of the so-called {\it dilation } method, see  \cite{Stein1},
\cite{Talenti1} alike the one-dimensional case \cite{Ostrovsky4}.\par

\vspace{3mm}

{\bf B. Weight Fourier transforms.}

\vspace{3mm}

 We consider in this subsection the multidimensional (vector): $l \ge 2 $ generalization of
PBO inequality. Namely, we investigate the inequality of a view

$$
| \ F_{\alpha,\beta}[f] (\vec{y}) \ |_{\vec{q}} \le K_{d, \vec{\alpha}, \vec{\beta} } (\vec{p})
|  \ f( \vec{x}) \ |_{\vec{p}} \eqno(2.2.1)
$$
or for simplicity

$$
| \ |y|^{-\alpha } \ F(y) \ |_q \le K_{d, \alpha, \beta} (p)
 \ | \ |x|^{\beta} \ f( x) \ |_{p}. \eqno(2.2.2)
$$

 Let us impose the following constrains:

$$
1 < p_j \le q_j < \infty, 0 \le \alpha_j < m_j/q_j, 0 \le \beta_j < m_j/p_j' \eqno(2.2.3)
$$
or equally $  p_j > m_j/(m_j-\beta_j); $
$$
[p' \stackrel{def}{=} p/(p-1) ], \   \beta_j - \alpha_j  = m_j(1 - 1/p_j - 1/q_j). \eqno(2.2.4)
$$

\vspace{3mm}

{\bf Theorem 2.2.} \par

{\bf 1.} The conditions (2.2.3) and (2.2.4) are necessary and sufficient for the existence and finiteness
of the constant $ K_{d, \vec{\alpha}, \vec{\beta} } (\vec{p}) $  for the inequality (2.2.1.) \par

{\bf 2.} If the conditions (2.2.3) and (2.2.4) are satisfied,
then the sharp (minimal) value of the coefficient $ K_{d, \vec{\alpha}, \vec{\beta} } (\vec{p}) $ satisfies the inequalities

$$
C_1(d, \vec{\alpha}, \vec{\beta}) \ \prod_{j=1}^l \left[\frac{p_j}{p_j -m_j/(m_j-\beta_j)} \right]^{(\alpha_j+\beta_j)/m_j  } \le
K_{d, \vec{\alpha}, \vec{\beta} } (\vec{p}) \le
$$

$$
C_2(d, \vec{\alpha}, \vec{\beta}) \
\prod_{j=1}^l \left[\frac{p_j}{p_j -m_j/(m_j-\beta_j)} \right]^{ \max(1,(\alpha_j+\beta_j)/m_j)  }. \eqno(2.2.5)
$$

{\bf Proof} of the second  proposition is at the same as the proof of of the second  proposition of theorem 2.1 and may
be omitted.  It based on the one-dimensional case $ l=1 $ and contained particularly in the articles \cite{Ostrovsky3},
\cite{Ostrovsky4}. \par

 The first assertion of theorem 2.2  may be proved again by the dilation method.   Indeed, we have in the case $ l = 1 $
and for $ \lambda = \const > 0: $

 $$
 |T_{\lambda} [f]|_p = \lambda^{-d/p} \ |f|_p; \  F_{\alpha,\beta} [T_{\lambda}f](x/\lambda) = \lambda^{-d + \alpha - \beta} \
 F_{\alpha, \beta} [f](x);
 $$

$$
F_{\alpha,\beta} [T_{\lambda}f](x) = \lambda^{-d + \alpha - \beta} \ F_{\alpha,\beta} [f](\lambda \ x);
$$

$$
| \ F_{\alpha, \beta} [T_{\lambda}f](\cdot) \ |_q = \lambda^{-d + \alpha - \beta - d/q} \  | \ F_{\alpha,\beta} [f](\cdot) \ |_q;;
$$
therefore

$$
-d + \alpha - \beta -d/q = -d/p, \ \alpha - \beta = d \left(1 + \frac{1}{q} - \frac{1}{p} \right).
$$

 In the multidimensional case we construct the "counter-example" functions as follows:

$$
f(\vec{x}) = \prod_{j=1}^l g_j(x_j), \  0 \ne g_j(\cdot) \in S(R^{m_j}),
$$
and obviously

$$
 \alpha_j - \beta_j = m_j \left(1 + \frac{1}{q_j} - \frac{1}{p_j} \right).
$$

\vspace{3mm}
{\bf Remark 2.2.1.} It is proved also in \cite{Ostrovsky3} that the classical {\it multidimensional} PBO
inequality, i.e. in the ordinary Lebesgue spaces $ L_p, $
see e.g. \cite{Okikiolu1}, p. 313-315; \cite{Gorbachev1} etc. may be obtained by virtue of equality $ |f|_{p,p} = |f|_p $
from theorem 2.2 as a particular case iff

$$
\frac{\beta_j - \alpha_j}{m_j} = \const, \ j = 1,2,\ldots,l.
$$

\vspace{3mm}

{\bf Remark 2.2.2.} The particular case $ \vec{\alpha}= \vec{\beta} =0 $ was considered by
Benedek A. and Panzone R. in the year 1961 \cite{Benedek1}; see also \cite{Gorbachev1}.\par

\vspace{3mm}

\hfill $\Box$ \\
 \bigskip

\section{ Multiple weight Riesz and Fourier transforms in anisotropic Grand Lebesgue-Riesz spaces.}

\vspace{3mm}

{\bf 0.}  Let $ Q $ be {\it appropriate} for concrete considered problem: Riesz potential (R) or
Fourier transform (F) convex (bounded or not) subset of the set $ \otimes_{j=1}^l [1,\infty]. $
 Let $ \psi = \psi(\vec{p}) $ be continuous in an interior $ Q^0 $ of the set $ Q $
strictly  positive  function such that

$$
\inf_{\vec{p} \in Q^0}  \psi(\vec{p}) > 0; \ \inf_{\vec{p} \notin Q^0}  \psi(\vec{p}) = \infty.
$$

 Let  $ f(\vec{x}) = f(x)  $ be some function from the space $ AG\psi. $  We denote for both the
considered problem  the introduced before correspondent function $  \vec{q} =\vec{q}(\vec{p})  $
and denote the inverse function by $  \vec{p} =\vec{p}(\vec{q}).  $ \\

\vspace{3mm}

{\bf Problem R. }  Define a new function

$$
\nu_R(\vec{q}) = \psi(\vec{p}(\vec{q})) \cdot K_{{\vec{m}; \vec{\alpha}, \vec{\beta}, \vec{\gamma}}}(\vec{p}(\vec{q})).
$$

{\bf Theorem 3.1.}

$$
|| I_{\vec{m};\vec{\alpha}, \vec{\beta}, \vec{\gamma}}[f](\vec{\cdot})||AG\nu_R \le
1 \cdot ||f||AG\psi, \eqno(3.1)
$$
where the constant $  "1" $ is the best possible.\\

\vspace{3mm}

{\bf Problem F. }  Define a new function

$$
\nu_F(\vec{q}) = \psi(\vec{p}(\vec{q})) \cdot K_{d; \vec{\alpha}, \vec{\beta}, \vec{\gamma}}(\vec{p}(\vec{q})).
$$

\vspace{3mm}

{\bf Theorem 3.2.}

$$
|| F_{\vec{m};\vec{\alpha}, \vec{\beta}}[f](\vec{\cdot})||AG\nu_F \le
1 \cdot ||f||AG\psi, \eqno(3.2)
$$
where the constant $  "1" $ is the best possible.\\

\vspace{3mm}

{\bf Proofs.}  It is sufficient to prove theorem 3.1; the second proposition (3.2) provided analogously. \par
Let  $  f(\cdot) \in AG\psi; $ we can suppose without loss of generality $ ||f||AG\psi = 1. $  This imply that

$$
|f|_{\vec{p}} \le \psi(\vec{p}).
$$
 We have using the result of theorem 2.1

 $$
 |u(\cdot)|_q \le K_{{\vec{m}; \vec{\alpha}, \vec{\beta}, \vec{\gamma}}}(\vec{p}) \cdot \psi( \vec{p})\le
 K_{{\vec{m}; \vec{\alpha}, \vec{\beta}, \vec{\gamma}}}(\vec{p}) \cdot \psi( \vec{p})\ ||f|| AG\psi. \eqno(3.3)
 $$
 As long as the variable $  \vec{p} $ is uniquely defined monotonic function on  $ \vec{q}, $  the inequality (3.2) is
equivalent to the assertion of theorem 3.1.\par
 The {\it exactness} of this estimation is proved in one-dimensional $ l=1 $ general case  in the article
\cite{Ostrovsky3};  the multidimensional case $ l \ge 2 $ provided analogously. \par

\hfill $\Box$ \\
 \bigskip

\section{ Multiple Riesz potential in bounded domains. }

\vspace{3mm}

We consider in this subsection the  {\it truncated}  Riesz's operator of a view

 $$
u^{(B)} = u^{(B)}(x)  = I^{(B)}_{\alpha,\beta,\gamma}[f](x) =I^{(B)}_{d;\alpha,\beta,\gamma}[f](x) = |x|^{-\beta}
\int_{B} \frac{f(y) \ |y|^{-\alpha} \ dy}{|x - y|^{\gamma}}, \eqno(4.1)
 $$
where $ B $ is open bounded domain in $ R^d $ contained the origin and such that

$$
0 < \inf_{x \in \partial B } |x| \le \sup_{x \in \partial B } |x| < \infty, \eqno(4.2)
$$
$ \partial B $ denotes boundary of the set $ B. $ \par

 For instance, we can  assume  that the set $ B $ is  unit ball in the space $ R^d $
with the center in origin and with finite positive radii:
 $$
B = \{x, \ x \in R^d, \ |x| < r \},  \ r = \const \in (0,\infty).
 $$

\vspace{3mm}

 {\bf Definition 4.1} We will call the subsets  $ \{ B \} $ satisfying formulated conditions (4.2)
 as {\it interior } domains. Notation: $ B \in I(R^d). $ \par
  In contradiction, the set $ D $ in the space $ R^d $ is said to be {\it exterior } domain if
by definition the open complement of $ D: \ (R^d \setminus D)^o $   is {\it interior } domain.
Notation: $ D \in E(R^d). $   \par

\vspace{3mm}

  It is proved in fact  in the article  \cite{Ostrovsky3} that in our notations and under our
assertions  if $ B \in I(R^d), $  then
there is a positive finite constant $  K^{(B)}_R(d; \alpha,\beta, \gamma) $ for which

$$
|I^{(B)}_{\alpha,\beta,\gamma}[f]|_q  \le K^{(B)}_{R;d; \alpha,\beta, \gamma}(p) \ |f|_p, \
p \in (1, p_+),  \eqno(4.3)
$$
where by definition

$$
K^{(B)}_{R;d; \alpha,\beta, \gamma}(p) =  \sup_{ p \in (1, p_+) } \sup_{0 < |f|_p   < \infty}
\left[\frac{|I^{(B)}_{\alpha,\beta,\gamma}[f]|_{q(p)}}{ |f|_p } \right].
$$

 It is known \cite{Ostrovsky3} that the coefficient $ K^{(B)}_R(d; \alpha,\beta, \gamma) $ is under formulated before
conditions bilateral bounded:

 $$
\frac{C(1)^{(B)}(R; d; \alpha,\beta, \gamma)}{\left[p_+ - p \right]^{\kappa}   }  \le K^{(B)}_{R;d; \alpha,\beta, \gamma}(p )
\le \frac{ C(2)^{(B)}(R;d; \alpha,\beta, \gamma)}{\left[p_+ - p \right]^{\kappa} }.  \eqno(4.4)
 $$

\vspace{3mm}

 We define as in section 2 the following {\it multiple} weight Riesz (linear) potential in bounded domain $  B: $:

 $$
u^{(B)}(\vec{x}) = I^{(B)}_{\vec{\alpha}, \vec{\beta}, \vec{\gamma}}[f](\vec{x})  =
 I^{(B)}_{\vec{m};\vec{\alpha}, \vec{\beta}, \vec{\gamma}}[f](\vec{x}) =
\otimes_{j=1}^l I^{(B)}_{m_j; \alpha_j,\beta_j,\gamma_j}[f](\vec{x})  =
$$

$$
 \int_{B}  \ \frac{|x_1|^{-\beta_1} \ |y_1|^{-\alpha_1} \ dy_1}{|x_1 - y_1|^{\gamma_1}}
\left[  \int_{B}  \ \frac{|x_2|^{-\beta_2} \ |y_2|^{-\alpha_2} \ dy_2}{|x_2 - y_2|^{\gamma_2}} \left[ \ldots
\left[ \int_{B}  \ \frac{|x_l|^{-\beta_l} \ |y_l|^{-\alpha_l} \ f( \vec{y} ) \ dy_l}{|x_l - y_l|^{\gamma_l}} \right]
 \right]   \right]. \eqno(4.5)
$$

  The following result based on the estimation (4.4) is obtained analogously to the proof of theorem 2.1. \\

\vspace{3mm}

{\bf Theorem 4.1a.} \par
 Let $ B \in I(R^d). $ The conditions (2.1.1) are necessary and sufficient for the existence of non-trivial coefficient
$ K^{(B)}_{{\vec{m}; \vec{\alpha}, \vec{\beta}, \vec{\gamma}}}(\vec{p}) $ for the following estimate:

$$
| I^{(B)}_{\vec{m};\vec{\alpha}, \vec{\beta}, \vec{\gamma}}[f](\vec{\cdot})|_{\vec{q}} \le
K^{(B)}_{{\vec{m}; \vec{\alpha}, \vec{\beta}, \vec{\gamma}}}(\vec{p}) \ |f(\cdot)|_{ \vec{p}}, \eqno(4.6)
$$
where for the minimal value of coefficient $  K^{(B)}_{{\vec{m}; \vec{\alpha}, \vec{\beta}, \vec{\gamma}}}(\vec{p}) $
holds the following inequality:

$$
 \frac{C^{(B)}_3(\vec{m}, \vec{\alpha}, \vec{\beta}, \vec{\gamma})}
 { \prod_{j=1}^l \left[(p_+^{(j)} -p_j) \right]^{\kappa_j} } \le
    K^{(B)}_{{\vec{m}; \vec{\alpha}, \vec{\beta}, \vec{\gamma}}}(\vec{p}) \le
 \frac{C^{(B)}_4(\vec{m}, \vec{\alpha}, \vec{\beta}, \vec{\gamma})}
 { \prod_{j=1}^l \left[(p_+^{(j)} -p_j) \right]^{\kappa_j} }, \eqno(4.7)
$$
if for all the values $ j \ \Rightarrow  p_j \in \left(1, p_+^{(j)} \right) $
 and $  K^{(B)}_{{\vec{m}; \vec{\alpha}, \vec{\beta}, \vec{\gamma}}}(\vec{p}) = \infty  $ in other case. \par

 For the exterior domain $ D $ we conclude: \\

{\bf Theorem 4.1b.} \par
 Let $ D \in E(R^d). $ The conditions (2.1.1) are necessary and sufficient for the existence of non-trivial coefficient
$ K^{(D)}_{{\vec{m}; \vec{\alpha}, \vec{\beta}, \vec{\gamma}}}(\vec{p}) $ for the following estimate:

$$
| I^{(D)}_{\vec{m};\vec{\alpha}, \vec{\beta}, \vec{\gamma}}[f](\vec{\cdot})|_{\vec{q}} \le
K^{(D)}_{{\vec{m}; \vec{\alpha}, \vec{\beta}, \vec{\gamma}}}(\vec{p}) \ |f(\cdot)|_{ \vec{p}}, \eqno(4.8)
$$
where for the minimal value of coefficient $  K^{(D)}_{{\vec{m}; \vec{\alpha}, \vec{\beta}, \vec{\gamma}}}(\vec{p}) $
holds the following inequality:

$$
 \frac{C^{(D)}_5(\vec{m}, \vec{\alpha}, \vec{\beta}, \vec{\gamma})}
 { \prod_{j=1}^l \left[(p_j -p_-^{(j)} \right]^{\kappa_j} } \le
    K^{(D)}_{{\vec{m}; \vec{\alpha}, \vec{\beta}, \vec{\gamma}}}(\vec{p}) \le
 \frac{C^{(D)}_6(\vec{m}, \vec{\alpha}, \vec{\beta}, \vec{\gamma})}
 { \prod_{j=1}^l \left[p_j -p_-^{(j)} \right]^{\kappa_j} }, \eqno(4.9)
$$
if for all the values $ j \ \Rightarrow  p_j \in \left( p_-^{(j)}, \infty \right) $
 and $  K^{(D)}_{{\vec{m}; \vec{\alpha}, \vec{\beta}, \vec{\gamma}}}(\vec{p}) = \infty  $ in other case. \par

\hfill $\Box$ \\
 \bigskip

\section{ Concluding remarks.}

\vspace{3mm}

{\bf A. Generalization on the regular varying weights.} \par

\vspace{3mm}

 We consider in this subsection some generalization on the {\it weight} Riesz's potential
 (and further on the Fourier weight operator)  of a view

 $$
 I_{d;\alpha, \delta}^{(S)} f(x) = \int_{R^d}
 \frac{f(y) \ | \log |x - y| \ |^{\delta} \ S(|\log|x - y| \ |) \ dy }
 { |x - y|^{d -\alpha} }, \eqno(5.1A.0)
 $$
where $ \alpha = \const \in (0,d), \ \delta = \const > 0, \ p \in (1, d/\alpha), $
$ p $ is the following function of variable $ q: $

$$
p = p(q) = \frac{dq}{d + \alpha q}, \eqno(5.1A.1)
$$
and conversely $ q = q(p), \ \vec{q} = \vec{q}(\vec{p});  $
$ S(z) $ is a {\it slowly varying} as $ z \to \infty $ continuous positive function:

$$
\forall x > 0 \ \Rightarrow  \lim_{z \to \infty} S(x z)/S(z) = 1.
$$

 We refer the reader to the book of Seneta \cite{Seneta1} to the using further facts about regular
and slowly varying functions.\par

 It is known, see \cite{Ostrovsky3} that under some simple conditions

$$
 |I_{d;\alpha,\delta}^{(S)} f|_q \le
 \frac{C(S) \ |f|_p } {[(p - 1) \ (d/\alpha - p)]^{1 + \delta - \alpha/d } }, \
 p \in (1, d/\alpha). \eqno(5.1A.2)
 $$
{\it and the last inequality is exact up to multiplicative constant.}\par

\vspace{4mm}

 We define as in section 2 the following {\it multiple}  Riesz (linear) potential with
regular varying weight

 $$
z^{(S)}(\vec{x}) = I^{(S)}_{\vec{\alpha}, \vec{\delta}}[f](\vec{x})  =
 I^{(S)}_{\vec{m};\vec{\alpha}, \vec{\delta}}[f](\vec{x}) =
\otimes_{j=1}^l I^{(S)}_{m_j; \alpha_j,\delta_j}[f](\vec{x}). \eqno(5.1A.3)
$$
where $ \alpha_j = \const \in (0,m_j), \ \delta_j = \const > 0, \ p_j \in (1, m_j/\alpha_j), $
$ p_j $ is the following function of variable $ q_j: $

$$
p_j = p_j(q_j) = \frac{m_j q_j}{m_j + \alpha_j q_j} \eqno(5.1A.4)
$$
and conversely $ q_j = q_j(p_j), \ \vec{q} = \vec{q}(\vec{p}).  $ \par
 $ \{ S_j(z) \} $ are  slowly varying as $ z \to \infty $ continuous positive functions.\par

  The following result based on the estimation  (5.1A.2) is obtained analogously the proof of theorem 2.1. \\

\vspace{3mm}

{\bf Theorem  5.1A.1. } \par
 The conditions  (5.1A.1) are necessary and sufficient for the existence non-trivial coefficient
$ K^{(S)}_{{\vec{m}; \vec{\alpha}, \vec{\delta}}}(\vec{p}) $ for the following estimate:

$$
| I^{(S)}_{\vec{m};\vec{\alpha}, \vec{\delta}}[f](\vec{\cdot})|_{\vec{q}} \le
K^{(S)}_{{\vec{m}; \vec{\alpha}, \vec{\delta}, }}(\vec{p}) \ |f(\cdot)|_{ \vec{p}}, \eqno(5.1A.4)
$$
where for the minimal value of coefficient $  K^{(S)}_{{\vec{m}; \vec{\alpha}, \vec{\delta}, }}(\vec{p}) $
holds the following inequality:

$$
 \frac{C^{(S)}_3(\vec{m}, \vec{\alpha}, \vec{\delta})}
 { \prod_{j=1}^l \left[ [(p_j - 1) \ (m_j/\alpha_j - p_j)]^{1 + \delta_j - \alpha_j/m_j } \right] }  \le
    K^{(S)}_{{\vec{m}; \vec{\alpha}, \vec{\delta}}}(\vec{p}) \le
$$

 $$
 \frac{C^{(S)}_4(\vec{m}, \vec{\alpha}, \vec{\delta})}
 { \prod_{j=1}^l \left[ [(p_j - 1) \ (m_j/\alpha_j - p_j)]^{1 + \delta_j - \alpha_j/m_j } \right] }.  \eqno(5.1A.5)
 $$
if for all the values $ j \ \Rightarrow  p_j \in \left(1, m_j/\alpha_j \right) $
 and $  K^{(S)}_{{\vec{m}; \vec{\alpha}, \vec{\delta}, }}(\vec{p}) = \infty  $ in other case. \par

\vspace{6mm}

 Analogously may be provided the multiple {\it regular varying weight}  generalization of PBO
inequality on the Anisotropic Lebesgue-Riesz spaces; ordinary case is investigated in \cite{Ostrovsky3}. \par

 Let $ L_j = L_j(z), M_j = M_j(z),  \ j = 1,2, \ldots, l, \ z \in (0,\infty) $ be a family of slowly varying
{\it simultaneously} as $ z \to 0 $ and as $ z \to \infty $ continuous positive functions:

$$
\lim_{z \to 0} \frac{L_j(x \ z)}{L_j(z) } = \lim_{z \to \infty} \frac{L_j(x \ z)}{L_j(z) } = 1;
$$

$$
\lim_{z \to 0} \frac{M_j(x \ z)}{M_j(z) } = \lim_{z \to \infty} \frac{M_j(x \ z)}{M_j(z) } = 1.
$$

Let us consider the following inequality:

$$
| \ | \vec{y}|^{-\vec{\alpha}} \  \prod_{j=1}^l L_j(|y_j|) \ F[f](\vec{y}) \ |_{\vec{q}} \le
K_{\vec{L}, \vec{M}, \vec{m}, \vec{\alpha}, \vec{\beta} } (\vec{p})
| \ |\vec{x}|^{\vec{\beta}} \ \prod_{j=1}^l M_j(|x_j|) \ f( \vec{x}) \ |_{\vec{p}}.  \eqno(5.2A.1)
$$

 Let us impose the following   restrictions:

$$
1 < p_j \le q_j < \infty, 0 \le \alpha_j < m_j/q_j, 0 \le \beta_j < m_j/p_j'; \eqno(5.2A.2)
$$
or equally $  p_j > m_j/(m_j-\beta_j); $
$$
  \beta_j - \alpha_j  = m_j(1 - 1/p_j - 1/q_j); \eqno(5.2A.3)
$$

$$
M_j(z) \asymp L_j(1/z) \Leftrightarrow \  0 < \inf_z \frac{M_j(z)}{L_j(1/z)}  \le \sup_z \frac{M_j(z)}{L_j(1/z)}  < \infty.
\eqno(5.2A.4)
$$

\vspace{3mm}

{\bf Theorem 5.2A.1.} \par

{\bf 1.} The conditions (5.2A.2), (5.2A.3) and (5.2A.4) are necessary and sufficient for the existence and finiteness
of the constant $ K_{d, \vec{\alpha}, \vec{\beta} } (\vec{p}) $  for the inequality (5.2A.1). \par
{\bf 2.} If the conditions  (5.2A.2), (5.2A.3) and (5.2A.4) are satisfied,
then the sharp (minimal) value of the coefficient $ K_{\vec{L}, \vec{M}, \vec{m}, \vec{\alpha}, \vec{\beta} } (\vec{p}) $
satisfies the inequalities

$$
C_1( \vec{L}, \vec{M}, \vec{m},\vec{\alpha}, \vec{\beta}) \
\prod_{j=1}^l \left[\frac{p_j}{p_j -m_j/(m_j-\beta_j)} \right]^{(\alpha_j+\beta_j)/m_j  } \le
K_{\vec{L}, \vec{M}, \vec{m}, \vec{\alpha}, \vec{\beta} } (\vec{p}) \le
$$

$$
C_2( \vec{L}, \vec{M}, \vec{m},\vec{\alpha}, \vec{\beta})
\prod_{j=1}^l \left[\frac{p_j}{p_j -m_j/(m_j-\beta_j)} \right]^{ \max(1,(\alpha_j+\beta_j)/m_j)  }. \eqno(5.2A.5)
$$

{\bf Proof} is at the same as the proof of theorem 2.1 and may be omitted.  It based on the one-dimensional case
$ l=1 $ and contained particularly in the articles \cite{Ostrovsky3},  \cite{Ostrovsky4}. \par

\vspace{4mm}

{\bf B.  Composed Riesz and Fourier operators.  }

\vspace{3mm}

  Let the set $ J := [1,2,\ldots,d]  $ consists in two non-empty {\it disjoint} subsets (partition):
 $ J = J(R) \cup J(F), \ J(R) \cap J(F) = \emptyset. $  We define the following {\it composed} linear Riesz-Fourier weight
 operator ( $ R $ = Riesz, $ F $ = Fourier) as follows:

 $$
W_{\vec{\alpha}, \vec{\beta}, \vec{\gamma}}[f](\vec{x})  \stackrel{def}{=} \left( \otimes_{j \in J(R) } I_{\alpha_j,\beta_j,\gamma_j} \right) \otimes
\left( \otimes_{j \in J(F)} F_{\alpha_j,\beta_j} \right) [f](\vec{x}). \eqno(5.1B.1)
 $$

 We investigate in this pilcrow  the inequality of a view

$$
| \ W_{\vec{\alpha}, \vec{\beta}, \vec{\gamma}}[f] \ |_{\vec{q}} \le K_{RF}(\vec{p}) \ | f|_{\vec{p}}. \eqno(5.1B.2)
$$

 We  impose the following conditions: \\

{\bf 1.} $ j \in J(F) \ \Rightarrow $

$$
1 < p_j \le q_j < \infty, 0 \le \alpha_j < m_j/q_j, 0 \le \beta_j < m_j/p_j'
$$
$$
  \beta_j - \alpha_j  = m_j(1 - 1/p_j - 1/q_j). \eqno(5.1B.3)
$$

{\bf 2.} $ j \in J(R) \ \Rightarrow $
 $$
 \gamma_j = \const, \alpha_j,\beta_j = \const \ge 0, \alpha_j + \gamma_j < m_j, \
 $$

$$
m_j+\frac{m_j}{q_j} = \frac{m_j}{p_j} + (\alpha_j+\beta_j+\gamma_j). \eqno(5.1B.4)
$$

  The relations (5.1B.3) and (5.1B.4) define uniquely determined functions $ \vec{p} = \vec{p}(\vec{q})   $
or conversely functions    $ \vec{q} = \vec{q}(\vec{p}).   $ \\

\vspace{3mm}

{\bf Proposition 5.1B.1}. \\

\vspace{3mm}

{\bf $ \alpha $}. The conditions (5.1B.3) and (5.1B.4) are necessary and sufficient for finiteness of the coefficient
$ K_{RF}(\vec{p})  $ in (5.1B.2).\\

{\bf $ \beta $}.  If this conditions are satisfied, then the minimal value of  coefficient $ K_{RF}(\vec{p})  $ is bilateral bounded
as follows:

$$
 C_5(\vec{m}, \vec{\alpha}, \vec{\beta}, \vec{\gamma})  \
\frac{ \prod_{j \in J(F)} ( p_j/( p_j - m_j/(m_j - \beta_j) ) )^{-(\alpha_j + \beta_j - \gamma_j)/m_j}}
{\prod_{j \in J(R)} [(p_j - p_-^{(j)}) \ (p_+^{(j)} - p_j )]^{ \kappa_j}} \le
$$

$$
K_{RF}(\vec{p}) \le C_6(\vec{m}, \vec{\alpha}, \vec{\beta}, \vec{\gamma})  \
\frac{ \prod_{j \in J(F)} ( p_j/( p_j -  m_j/(m_j - \beta_j) ) )^{-\max(1, (\alpha_j + \beta_j - \gamma_j)/m_j)}}
{\prod_{j \in J(R)} [(p_j - p_-^{(j)}) \ (p_+^{(j)} - p_j )]^{ \kappa_j}}.  \eqno(5.1B.5)
$$

 This assertion  may be simple obtained by synthesis of theorems 2.1 and 2.2. \par

\vspace{3mm}

{\bf C.  Mixture Riesz and Fourier operators.  }

\vspace{3mm}

  The following (linear) operator

$$
G_{\alpha,\beta, \gamma}[f](x) = |x|^{-\beta} \ \int_{R^d} e^{i xy} \ \frac{|y|^{-\alpha} \ f(y) \ dy}{ |x-y|^{\gamma}}, \
 x,y \in R^d \eqno(5.1C.1)
$$
is said to be "one-dimensional",  i.e.  $  l=1  $ mixture Riesz and Fourier operator. \par
 We consider also in this subsection its multidimensional $ \ l \ge 2 $ generalization

 $$
v(\vec{x}) = G_{\vec{\alpha}, \vec{\beta}, \vec{\gamma}}[f](\vec{x})  =
 G_{\vec{m};\vec{\alpha}, \vec{\beta}, \vec{\gamma}}[f](\vec{x}) =
\otimes_{j=1}^l G_{m_j; \alpha_j,\beta_j,\gamma_j}[f](\vec{x})  =
$$

$$
 \int_{R^{m_1}}  \  e^{i x_1 y_1} \ \frac{|x_1|^{-\beta_1} \ |y_1|^{-\alpha_1} \ dy_1}{|x_1 - y_1|^{\gamma_1}}
[ \ \int_{R^{m_2}} \  e^{i x_2 y_2}  \ \frac{|x_2|^{-\beta_2} \ |y_2|^{-\alpha_2} \ dy_2}{|x_2 - y_2|^{\gamma_2}} [ \ldots
$$

$$
[ \ \int_{R^{m_l}} \  e^{i x_l y_l}  \ \frac{|x_l|^{-\beta_l} \ |y_l|^{-\alpha_l} \ f( \vec{y} ) \ dy_l}{|x_l - y_l|^{\gamma_l}} \ ]
] \  ].  \eqno(5.1C.2)
$$

 The inequality of a view (relative the Lebesgue measures in whole space $ R^d $ )

$$
|G_{d;\alpha,\beta,\gamma}[f]|_q \le K_{GPBO}(p) \ |f|_p  \eqno(5.1C.3)
$$
is said to be generalized (G) weight Pitt - Beckner - Okikiolu (GPBO) inequality.\par

 We will understood as customary as the variable $ K_{GPBO}(p) $ its minimal value:

$$
K_{GPBO}(p) = \sup_{f \ne 0, |f|_p < \infty} \left[ \frac{|G_{d;\alpha,\beta,\gamma}[f]|_q}{|f|_p} \right].
$$

 Denote

 $$
 \tilde{p}_- = \frac{d}{d-\beta}, \ \tilde{q}_+ = \frac{d}{\alpha-\gamma},
 \tilde{p}_+ =\infty, \ \tilde{q}_- = 1, \eqno(5.1C.4)
 $$
 and suppose

$$
\beta + \gamma - \alpha = d \left(1-\frac{1}{p} - \frac{1}{q}  \right),
$$

 $$
 \alpha - \gamma > 0, \ \gamma < d, \ \alpha + \beta > \gamma,
 $$

 $$
 p \in  (\tilde{p}_-,  \tilde{p}_+) \ (\Leftrightarrow  q \in  (\tilde{q}_-,  \tilde{q}_+)), \
  p  < q.  \eqno(5.1C.5)
 $$

\vspace{3mm}

{\bf Proposition 5.1C.1}.

\vspace{3mm}

$$
C_1(d,\alpha,\beta, \gamma)  \left[\frac{p}{p-\tilde{p}_-} \right]^{(\alpha +\beta - \gamma)/d } \le K_{GPBO}(p) \le
$$

$$
C_2(d,\alpha,\beta,\gamma)  \left[\frac{p}{p-\tilde{p}_-} \right]^{ \max(1,(\alpha +\beta - \gamma)/d) }, \
 0 < C_{1,2}(d;\alpha,\beta,\gamma) < \infty.   \eqno(5.1C.6)
$$

 The right hand side of this bilateral estimate  may be obtained after  some computations
from the articles \cite{Jurkat1},  \cite{Schechter1}, \cite{Lacey1}, \cite{Sawyer1}, \cite{Wheeden1}.
  The opposite estimate may be obtained by means of at the same counterexample as
 in article \cite{Ostrovsky3}. \par
 Note that both the possibilities in the power in the right hand side in  (5.1C.6): "1" and
$ "(\alpha +\beta - \gamma)/d" $  are attainable. \par

\vspace{3mm}

{\it We consider in this subsection the multidimensional (vector): $l \ge 2 $ generalization of
GPBO inequality.} \par

 Namely, we investigate the inequality of a view

$$
| G_{\vec{m},\vec{\alpha}, \vec{\beta}, \vec{\gamma}}[f](\cdot)|_{\vec{q}} \le
K_{\vec{m}, \vec{\alpha}, \vec{\beta},\vec{\gamma} } (\vec{p}) \ |f|_{\vec{p}},
$$
where as usually

$$
K_{\vec{m}, \vec{\alpha}, \vec{\beta},\vec{\gamma} } (\vec{p}) \stackrel{def}{=} \sup_{0 < |f|_p < \infty}
\left[\frac{| G_{\vec{m},\vec{\alpha}, \vec{\beta}, \vec{\gamma}}[f](\cdot)|_{\vec{q}}}{|f|_{\vec{p}}} \right].
\eqno(5.1C.7)
$$

  Denote as before

 $$
 \tilde{p}_-^{(j)} = \frac{m_j}{m_j-\beta_j}, \ \tilde{q}_+^{(j)} = \frac{m_j}{\alpha_j-\gamma_j},
 \tilde{p}_+^{(j)} =\infty, \ \tilde{q}_-^{(j)} = 1, \eqno(5.1C.8)
 $$
 and suppose

$$
\beta_j + \gamma_j - \alpha_j = m_j \left(1-\frac{1}{p_j} - \frac{1}{q_j}  \right),
$$

 $$
 \alpha_j - \gamma_j > 0, \ \gamma_j < m_j, \ \alpha_j + \beta_j > \gamma_j,
 $$

 $$
 p_j \in  (\tilde{p}_-^{(j)},  \tilde{p}_+^{(j)} ) \ (\Leftrightarrow  q_j \in  (\tilde{q}_-^{(j)},  \tilde{q}_+^{(j)})), \
  p_j  \le q_j,  \eqno(5.1C.9)
 $$
so that in equality (5.1C.7) the vector $ \vec{q} $ is uniquely defined function on the vector $ \vec{p}: \
 \vec{q} = \vec{q}(\vec{p}). $ \par

 \vspace{3mm}

{\bf Theorem 5.1C.1.} \\
{\bf 1.} The conditions  (5.1C.9) are necessary and sufficient for finiteness of the coefficient
$ K_{\vec{m}, \vec{\alpha}, \vec{\beta},\vec{\gamma} } (\vec{p}). $\\
{\bf 2.} If the conditions  (5.1C.9) are satisfied, then

$$
C_1(\vec{m}, \vec{\alpha}, \vec{\beta},\vec{\gamma} ) \prod_{j=1}^l \left[ \frac{p_j}{p_j - \tilde{p}_-^{(j)}} \right]^{(\alpha_j + \beta_j - \gamma_j)/m_j } \le K_{\vec{m}, \vec{\alpha}, \vec{\beta},\vec{\gamma} } (\vec{p}) \le
$$

$$
C_2(\vec{m}, \vec{\alpha}, \vec{\beta},\vec{\gamma} )
\prod_{j=1}^l \left[ \frac{p_j}{p_j - \tilde{p}_-^{(j)}} \right]^{\max(1,(\alpha_j + \beta_j - \gamma_j)/m_j) }. \eqno(5.1C.10)
$$

 \vspace{3mm}

 {\bf Sketch of proof.} The necessity of the conditions  (5.1C.9)  may be grounded as before  by means of
multidimensional dilation method; see e.g. \cite{Ostrovsky3}. The estimates (5.1C.10) may be proved  as
in the theorem 2.1 on the basis the one-dimensional version. \par

\hfill $\Box$ \\

\end{document}